\documentclass[12pt]{article}%
\usepackage{amsfonts}
\usepackage{amsmath,amssymb}
\usepackage[applemac]{inputenc}
\usepackage{amsmath,amssymb,fullpage}
\usepackage{showlabels}
\usepackage{color}
\usepackage{amsmath}
\usepackage{amssymb}
\usepackage{graphicx}%
\providecommand{\U}[1]{\protect\rule{.1in}{.1in}}

\newtheorem{definition}{Definition}[section]
\newtheorem{example}{Example}[section]
\newtheorem{theorem}[definition]{Theorem}
\newtheorem{problem}[definition]{Problem}
\newtheorem{remark}[definition]{ \it Remark}

\newtheorem{lemma}[definition]{Lemma}

\numberwithin{equation}{section}

\def\1B{\text{1\!\!I}}

\begin{document}

\title{Mean-Field Stochastic Control with Elephant Memory in Finite and Infinite Time Horizon}
\author{Nacira AGRAM$^{1,2}$ and Bernt ØKSENDAL$^{1}$}
\date{21 June 2019
\vskip 0.2cm
The final version of this paper will be published in Stochastics}
\maketitle

\begin{abstract}
\noindent Our purpose of this paper is to study stochastic control problems
for systems driven by mean-field stochastic differential equations with
elephant memory, in the sense that the system (like the elephants) never
forgets its history. We study both the finite horizon case and the infinite
time horizon case.

\begin{itemize}
\item In the finite horizon case, results about existence and uniqueness of
solutions of such a system are given. Moreover, we prove sufficient as well as
necessary stochastic maximum principles for the optimal control of such
systems. We apply our results to solve a mean-field linear quadratic control problem.

\item For infinite horizon, we derive sufficient and necessary maximum
principles.\newline As an illustration, we solve an optimal consumption
problem from a cash flow modelled by an elephant memory mean-field system.

\end{itemize}
\end{abstract}

\paragraph{MSC(2010):}

60H05, 60H20, 60J75, 93E20, 91G80,91B70.

\paragraph{Keywords:}

Mean-field stochastic differential equation; memory; stochastic maximum
principle; partial information; backward stochastic differential equation.

\footnotetext[1]{Department of Mathematics, University of Oslo, P.O. Box 1053
Blindern, N--0316 Oslo, Norway. Email: \texttt{naciraa@math.uio.no,
oksendal@math.uio.no.}
\par
This research was carried out with support of the Norwegian Research Council,
within the research project Challenges in Stochastic Control, Information and
Applications (STOCONINF), project number 250768/F20.}

\footnotetext[2]{University Mohamed Khider of Biskra, Algeria.}

\section{Introduction}

In this paper we study optimal control of stochastic systems with memory.
There are many ways of modelling such systems. Examples include systems with
delay or Volterra integral equations. See e.g. Agram \textit{et al}
\cite{AOY}. Here we are interested in stochastic differential equations (SDEs)
where the coefficients of the system depend upon the whole past. In this case
we say that the system has \emph{elephant memory}, inspired by the folklore
that an elephant never forgets. In addition we allow the dynamics of the
system to depend on both the current and previous laws of the state.
Specifically, we assume that the state $X(t)$ at time $t$ satisfies the
following equation%

\begin{equation}
\left\{
\begin{array}
[c]{ll}%
dX(t) & =b(t,X(t),X_{t},M(t),M_{t})dt+\sigma(t,X(t),X_{t},M(t),M_{t})dB(t)\\
& \text{ \ \ \ \ }+%
{\textstyle\int_{\mathbb{R}_{0}}}
\gamma(t,X(t),X_{t},M(t),M_{t},\zeta)\tilde{N}(dt,d\zeta);t\geq0,\\
X(0) & =x_{0},
\end{array}
\right.  \label{eq1.1}%
\end{equation}
where $X_{t}:=\{X(t-s)\}_{0\leq s\leq t}$ is the path of $X$ up to time $t$,
$M(t)=\mathcal{L}(X(t))$ is the law of $X(t)$, and $M_{t}:=\{M(t-s)\}_{0\leq
s\leq t}$ is the path of the law process.\newline We call equation
\eqref{eq1.1} a mean-field stochastic differential equation (MF-SDE) with
elephant memory. For more information on mean-field SDEs without memory we
refer to e.g. Carmona and Delarue \cite{CD1},\cite{CD2} and the references
therein. 

A historical process $X_{t}:=\{X(s)\}_{0\leq s\leq t}$ was studied by Dynkin
\cite{D}, but in a different framework. Different types of systems with
memories were discussed in the seminal work of Mohammed \cite{M}. A stochastic
version of Pontryagin's maximum principle for systems with delay
(discrete/distributed) has been derived by Chen and Wu \cite{CW}, Dahl
\textit{et al} \cite{DMOR} and Øksendal \textit{et al} \cite{osz}.\ 

The above mentioned works deal only with the finite horizon case. We refer to
Agram \textit{et al} \cite{AHOP}, \cite{AO} for the infinite time horizon setting.

Systems with discrete delay and mean-field have been studied by Meng and Shen
\cite{MS}, Agram and Røse \cite{AR}, but the mean-field terms considered there
are of a special kind, nameley the expectation of a function of the state,
i.e. $\mathbb{E}[\varphi(X(t-\delta))]$ for some bounded function $\varphi$
and $\delta$ is a positive delay constant.

In this paper we consider a more general situation, where the dynamics of the
state $X(t)$ at time $t$ depends on both the history of the state, the law for
the random variable $X(t)$ and the history of this law, as we have seen in
\eqref{eq1.1}. Moreover, we consider both the finite horizon case (Section 3)
and the infinite horizon case (Section 4).

Since the system is not Markovian, it is not obvious how to derive the dynamic
programming approach, but one can still get the HJB equation by using the
minimal backward stochastic differential equation (BSDE). This has been
studied by Fuhrman and Pham \cite{FP} by using the control randomization
method, considering measures defined on the Wasserstein metric space of
probability measures with finite second moment and using Lions lifting
techniques for differentiating the function of the measure.

In our paper, we use the Hilbert space of measures constructed in Agram
\textit{et al} \cite{AO3}, \cite{AO2}, \cite{ABOP}.

In Section 3 we obtain finite horizon maximum principles for the optimal
control of such systems. This is related to the paper by Agram and Øksendal in
\cite{AO2}, where the memorized paths are defined as $\{X(s)\}_{s\in\lbrack
t-\delta,t]}$ for a fixed $\delta>0$. However, in the current paper, we
consider as memory the whole trajectory $\{X(s)\}_{s\in\lbrack0,t]}.$

In the infinite horizon case in Section 4, we show that by replacing the
terminal value of the BSDE for the adjoint processes with a suitable
transversality condition at infinity, we can derive stochastic maximum
principles also in this case. As an illustration we study an infinite horizon
version of an optimal consumption problem with elephant memory. \color{black}

\section{Framework}

We now explain our setup in more detail:\newline Let $B=(B(t))_{t\in
\lbrack0,T]}$ and $\tilde{N}(dt,d\zeta)$ be a $d$-dimensional Brownian motion
and a compensated Poisson random measure, respectively, defined in a complete
filtered probability space $(\Omega,\mathcal{F},\mathbb{F},\mathbb{P}).$ The
filtration $\mathbb{F=}\left\{  \mathcal{F}_{t}\right\}  _{t\geq0}$ is assumed
to be the $\mathbb{P}$-augmented filtration generated by $B$ and $\tilde{N}$.
\newline

\subsection{Sobolev spaces of measures}

We now define a \emph{weighted Sobolev spaces of measures}. It is strongly
related to the space introduced in Agram and Øksendal \cite{AO3}, \cite{AO2},
but with a different weight, which is more suitable for estimates (see e.g.
Lemma 2.4 below):

\begin{itemize}
\item Let $n$ be a given integer. Then we define $\tilde{\mathcal{M}}%
=\tilde{\mathcal{M}}^{n}$ to be the pre-Hilbert space of random measures $\mu$
on $\mathbb{R}$ equipped with the norm
\[%
\begin{array}
[c]{lll}%
\left\Vert \mu\right\Vert _{\tilde{\mathcal{M}}^{n}}^{2} & := & \mathbb{E[}%
{\textstyle\int_{\mathbb{R}}}|\hat{\mu}(y)|^{2}(1+|y|)^{-n}dy]\text{,}%
\end{array}
\]
where $\hat{\mu}$ is the Fourier transform of the measure $\mu$, i.e.%
\[%
\begin{array}
[c]{lll}%
\hat{\mu}(y) & := & {\textstyle\int_{\mathbb{R}}}\color{black}e^{-ixy}%
d\mu(x);\quad y\in\mathbb{R}.
\end{array}
\]

\item For simplicity of notation, we will in the following fix
\[
n\geq2
\]
and we let $\mathcal{M}=\mathcal{M}^{n}$ denote the completion of
$\tilde{\mathcal{M}}=\tilde{\mathcal{M}}^{n}$ and we let $\mathcal{M}_{0}$
denote the set of deterministic elements of $\mathcal{M}$. 

\item Let $\mathcal{\tilde{M}}_{t}$ be the pre-Hilbert space of all paths
$\bar{\mu}=\{\mu(s)\}_{s\in\lbrack0,t]}$ of processes $\mu(\cdot)$ with
$\mu(s)\in\mathcal{\tilde{M}}^{n}\mathcal{=\tilde{M}}$ for each $s\in
\lbrack0,t]$ equipped with the norm%
\[
\left\Vert \bar{\mu}\right\Vert _{\mathcal{\tilde{M}}_{t}}^{2}:=%
{\textstyle\int_{0}^{t}}
||\mu(s)||_{\mathcal{\tilde{M}}}^{2}ds.
\]

\item We denote by $\mathcal{\tilde{M}}_{0,t}$ the set of all deterministic
elements of $\mathcal{\tilde{M}}_{t}$ and by $\mathcal{M}_{t}$ and
$\mathcal{M}_{0,t}$ their completions respectively.

\item If $\bar{x}\in\mathbb{R}^{[0,\infty)}$ (the set of all functions from
$[0,\infty)$ into $\mathbb{R}$), we define $\bar{x}_{t}\in\mathbb{R}%
^{[0,\infty)}$ by
\begin{align*}
&  \bar{x}_{t}(s)=\bar{x}(t-s);\quad s\in\lbrack0,t],\\
&  \bar{x}_{t}(s)=0;\quad s>t.
\end{align*}

\item If $\bar{x}\in\mathbb{R}^{[0,\infty)}$, we define $\bar{x}^{t}%
\in\mathbb{R}^{[0,\infty)}$ by
\begin{align}
&  \bar{x}^{t}(s)=\bar{x}(t+s);\quad s\in\lbrack0,t],\nonumber\\
&  \bar{x}^{t}(s)=0;\quad s>t. \label{fs}%
\end{align}

\end{itemize}

\noindent The following results is essential for our approach:

\begin{lemma}
\label{Lemma 4} Assume that $n \geq4$.

\begin{description}
\item[(i)] Let $X^{(1)}$ and $X^{(2)}$ be two $1$-dimensional random variables
in $L^{2}(\mathbb{P})$. \newline Then there exists a constant $C_{0}$ not
depending on $X^{(1)}$ and $X^{(2)}$ such that
\[%
\begin{array}
[c]{lll}%
\left\Vert \mathcal{L}(X^{(1)})-\mathcal{L}(X^{(2)})\right\Vert _{\mathcal{M}%
_{0}}^{2} & \leq & C_{0}\ \mathbb{E}[(X^{(1)}-X^{(2)})^{2}]\text{.}%
\end{array}
\]

\item[(ii)] Let $\{X^{(1)}(t)\}_{t\geq0},$ $\{X^{(2)}(t)\}_{t\geq0}$ be two
processes such that
\[
\mathbb{E}[%
{\textstyle\int_{0}^{T}}
X^{(i)2}(s)ds]<\infty\text{ for } i=1,2.
\]
Then, for all $t$,
\[%
\begin{array}
[c]{lll}%
||\mathcal{L}(X_{t}^{(1)})-\mathcal{L}(X_{t}^{(2)})||_{\mathcal{M}_{0,t}}^{2}
& \leq & C_{0}\ \mathbb{E}[%
{\textstyle\int_{0}^{t}}
(X^{(1)}(t-s)-X^{(2)}(t-s))^{2}ds]\text{.}%
\end{array}
\]

\end{description}
\end{lemma}

\noindent{Proof.} \quad By definition of the norms and standard properties of
the complex exponential function, we have%
\[%
\begin{array}
[c]{l}%
\left\Vert \mathcal{L}(X^{(1)})-\mathcal{L}(X^{(2)})\right\Vert _{\mathcal{M}%
_{0}}^{2}\\
={\textstyle\int_{\mathbb{R}}}|\widehat{\mathcal{L}}(X^{(1)})(y)-\widehat
{\mathcal{L}}(X^{(2)})(y)|^{2}(1+|y|)^{-n}dy\\
={\textstyle\int_{\mathbb{R}}}|{\textstyle\int_{\mathbb{R}}}e^{-ixy}%
d\mathcal{L}(X^{(1)})(x)-{\textstyle\int_{\mathbb{R}^{d}}}e^{-ixy}%
d\mathcal{L}(X^{(2)})(x)|^{2}(1+|y|)^{-n}dy\\
={\textstyle\int_{\mathbb{R}}}|\mathbb{E}[e^{-iyX^{(1)}}-e^{-iyX^{(2)}}%
]|^{2}(1+|y|)^{-n}dy\\
\leq{\textstyle\int_{\mathbb{R}}}\mathbb{E}[|e^{-iyX^{(1)}}-e^{-iyX^{(2)}%
}|^{2}](1+|y|)^{-n}dy\\
\leq{\textstyle\int_{\mathbb{R}}}y^{2}(1+|y|)^{-n}dy\mathbb{E}[|X^{(1)}%
-X^{(2)}|]^{2}\\
\leq C_{0}\mathbb{E}[(X^{(1)}-X^{(2)})^{2}],
\end{array}
\]
where
\begin{equation}
C_{0}=%
{\textstyle\int_{\mathbb{R}}}
y^{2}(1+|y])^{-n}dy<\infty\text{ since }n\geq4.
\end{equation}
Similarly we get that
\[%
\begin{array}
[c]{lll}%
||\mathcal{L}(X_{t}^{(1)})-\mathcal{L}(X_{t}^{(2)})||_{\mathcal{M}_{0,t}}^{2}
& \leq &
{\textstyle\int_{0}^{t}}
\left\Vert \mathcal{L}(X^{(1)}(t-s))-\mathcal{L}(X^{(2)}(t-s))\right\Vert
_{\mathcal{M}_{0}}^{2}ds\\
& \leq & C_{0}\ \mathbb{E}[%
{\textstyle\int_{0}^{t}}
(X^{(1)}(t-s)-X^{(2)}(t-s))^{2}ds].
\end{array}
\]
\hfill$\square$

\begin{remark}
If $f\in L^{1}(\mathbb{R})\cap L^{2}(\mathbb{R})$, then by the Fourier
inversion formula and the Fubini theorem
\begin{align}
\left\vert
{\textstyle\int_{\mathbb{R}}}
f(x)d\mu(x)\right\vert ^{2} &  =\left\vert
{\displaystyle\int_{\mathbb{R}}}
\left(
{\displaystyle\int_{\mathbb{R}}}
\frac{1}{2\pi}\widehat{f}(y)e^{ixy}dy\right)  d\mu(x)\right\vert
^{2}=\left\vert
{\displaystyle\int_{\mathbb{R}}}
\left(
{\textstyle\int_{\mathbb{R}}}
\frac{1}{2\pi}e^{ixy}d\mu(x))\widehat{f}(y)dy\right)  \right\vert
^{2}\nonumber\label{eq1.2}\\
&  =\frac{1}{2\pi}\left\vert
{\displaystyle\int_{\mathbb{R}}}
\widehat{\mu}(-y)\widehat{f}(y)dy\right\vert ^{2}.
\end{align}

If we assume in addition that $f^{\prime}\in L^{1}(\mathbb{R})\cap
L^{2}(\mathbb{R})$ then we know that
\begin{equation}%
{\displaystyle\int_{\mathbb{R}}}
\left\vert \widehat{f}(y)\right\vert ^{2}(1+|y|)^{2}dy\leq C(||\widehat
{f}||_{L^{2}(\mathbb{R})}^{2}+||\widehat{f^{\prime}}||_{L^{2}(\mathbb{R})}%
^{2})=:C_{2}(f)<\infty.
\end{equation}
and hence by \eqref{eq1.2} we get the following result:

\begin{lemma}
Suppose $f,f^{\prime}\in L^{1}(\mathbb{R})\cap L^{2}(\mathbb{R})$.Then
\[
|E[f(X(t))]|^{2}\leq\frac{1}{2\pi}C_{2}(f)||\mu||_{\mathcal{M}_{0}^{2}}^{2}.
\]

\end{lemma}
\end{remark}

\noindent{Proof.} \quad\
\begin{align*}
|E[f(X(t))]|^{2} &  =\left\vert
{\displaystyle\int_{\mathbb{R}}}
f(x)d\mu(x)\right\vert ^{2}\leq\frac{1}{2\pi}\left\vert
{\displaystyle\int_{\mathbb{R}}}
\widehat{\mu}(-y)\widehat{f}(y)dy\right\vert ^{2}\\
&  \leq\frac{1}{2\pi}{\textstyle%
{\displaystyle\int_{\mathbb{R}}}
}|\widehat{\mu}(y)|^{2}(1+|y|)^{-2}dy%
{\displaystyle\int_{\mathbb{R}}}
|\widehat{f}(y)|^{2}(1+|y|)^{2}dy\\
&  =\frac{1}{2\pi}C_{2}(f)||\mu||_{\mathcal{M}_{0}^{2}}^{2}.
\end{align*}

$\square$

This is a useful estimate, because if $\mu:=\mathcal{L}_{X(t)}$ where $X(t)$
solves a MF-SDE of the type \eqref{eq1.1}, then we always have $||\mu
||_{\mathcal{M}_{0}^{2}}<\infty$.\newline

Applying the previous result to $\mu:=\mu_{1}-\mu_{2}$ where $\mu
_{i}=\mathcal{L}(X_{i}(t));i=1,2,$ we the following Lipschitz estimate. This
is useful when we want to verify the Lipschitz condition (ii) in Section 3.1
in specific MF-SDEs with memory.

\begin{lemma}
Let $X_{1}(t),X_{2}(t)$ be two solutions of a MF-SDE, with corresponding laws
$\mu_{1},\mu_{2}$ at time $t$. Then if $f,f^{\prime}\in L^{1}(\mathbb{R})\cap
L^{2}(\mathbb{R})$, the following Lipschitz continuity holds:
\begin{equation}
|E[f(X_{1}(t))]-E[f(X_{2}(t))]|^{2}\leq\frac{1}{2\pi}C(f)||\mu_{1}-\mu
_{2}||_{\mathcal{M}_{0}^{2}}^{2}.
\end{equation}

\end{lemma}

\begin{definition}
[Law process]From now on we use the notation
\[
M_{t}:=M(t):=\mathcal{L}(X(t));\quad0\leq t\leq T
\]
for the law process $\mathcal{L}(X(t))$ of $X(t)$ with respect to $\mathbb{P}$.
\end{definition}

We recall the following result, which is proved in Agram and Øksendal
\cite{AO3}, Lemma 2.3/Lemma 6:

\begin{lemma}
\label{2.6} If $X(t)$ is an Itô-Lévy process as in \eqref{eq1.1}, then the map
$t\mapsto M(t):[0,T]\rightarrow\mathcal{M}_{0}^{4}$ is absolutely continuous.
Hence the derivative%

\[
M^{\prime}(t):=\frac{d}{dt}M(t)
\]
exists for a.a. $t$. and we have
\[
M(t)=M(0)+%
{\textstyle\int_{0}^{t}}
M^{\prime}(s)ds;\quad t\geq0.
\]

\end{lemma}

We will also use the following spaces:

\begin{itemize}
\item $C^{d}$ stands for the space of $\mathbb{R}^{d}$-valued continuous
functions defined over the time interval $[0,T].$

\item Given a finite time horizon $T>0$, for $1\leq p<+\infty,$ let
$S^{p}[0,T]$ denote the space of $%
\mathbb{R}
^{d}$-valued $\mathbb{F}$-adapted càdlàg processes $X=(X(t))_{t\in\left[
0,T\right]  }$ such that:%
\[
||X||_{S^{p}[0,T]}^{p}:=\mathbb{E[}\underset{t\in\left[  0,T\right]
}{\text{sup}}|X(t)|^{p}]<\infty.
\]

\item We define $\mathcal{\bar{S}}[0,T]$ the space of processes $\bar
{x}=\{x(s)\}_{0\leq s\leq t}:[0,T]\mapsto\mathbb{R}$ such that
\[
||\bar{x}||_{\mathcal{\bar{S}}[0,T]}^{2}:=\mathbb{E}[\underset{s\in\left[
0,t\right]  }{\text{sup}}x^{2}(s)ds]<\infty.
\]
For finite $T$ we identify functions $\bar{x}:[0,T]\mapsto\mathbb{R}$ with
functions $\bar{x}\in\mathcal{\bar{S}}[0,T]$ such that $x(s)=0$ for $s>T$, and
we regard them as functions defined on all $(-\infty,\infty)$ by setting
$x(s)=0$ for $s<0$.

\item Let $\mathcal{\bar{S}}[0,\infty)$ denote the space of processes $\bar
{x}=\{x(s)\}_{0\leq s\leq\infty}:[0,\infty)\mapsto\mathbb{R}$ such that
\[
||\bar{x}||_{\mathcal{\bar{S}}[0,\infty)}^{2}:=\mathbb{E}[{%
{\textstyle\int_{0}^{\infty}}
}x^{2}(s)ds]<\infty.
\]

\item We let $\mathbb{G}:=\{\mathcal{G}_{t}\}_{t\geq0}$ be a fixed given
subfiltration of $\mathbb{F}$ with $\mathcal{G}_{t}\subseteq\mathcal{F}_{t}$
for all $t\geq0$.The sigma-algebra $\mathcal{G}_{t}$ represents the
information available to the controller at time $t$. By $U$ we denote a
nonempty convex subset of $\mathbb{R}^{d}$ and we denote by $\mathcal{U}%
_{adm}$ the set of paths $U$-valued $\mathbb{G}$-predictable control
processes. We consider them as the \emph{admissible} control processes. \color{black}
\end{itemize}

\subsection{Fr\'echet derivatives and dual operators}

\noindent In this subsection we review briefly the Fréchet differentiability
and we introduce some dual operators, which will be used when we in the next
sections study Pontryagin's maximal principles for our stochastic control
problem.\newline Let $\mathcal{X},\mathcal{Y}$ be two Banach spaces and
let $F:\mathcal{X}\rightarrow\mathcal{Y}$. Then

\begin{itemize}
\item We say that $F$ has a directional derivative (or Gâteaux derivative) at
$v\in\mathcal{X}$ in the direction $w\in\mathcal{X}$ if
\[
D_{w}F(v):=\lim_{\varepsilon\rightarrow0}\frac{1}{\varepsilon}(F(v+\varepsilon
w)-F(v))
\]
exists in $\mathcal{Y}$.

\item We say that $F$ is Fréchet differentiable at $v\in\mathcal{X}$ if there
exists a continuous linear map $A:\mathcal{X}\rightarrow\mathcal{Y}$ such
that
\[
\lim_{\substack{h\rightarrow0\\h\in\mathcal{X}}}\frac{1}{\Vert h\Vert
_{\mathcal{X}}}\Vert F(v+h)-F(v)-A(h)\Vert_{\mathcal{Y}}=0,
\]
where $A(h)=\langle A,h\rangle$ is the action of the liner operator $A$ on
$h$. In this case we call $A$ the \textit{gradient} (or Fréchet derivative) of
$F$ at $v$ and we write
\[
A=\nabla_{v}F.
\]

\item If $F$ is Fréchet differentiable at $v$ with Fréchet derivative
$\nabla_{v}F$, then $F$ has a directional derivative in all directions
$w\in\mathcal{X}$ and
\[
D_{w}F(v)=\nabla_{v}F(w)=\langle\nabla_{v}F,w\rangle.
\]

\end{itemize}

\noindent In particular, note that if $F$ is a linear operator, then
$\nabla_{v}F=F$ for all $v$. \vskip0.2cm \emph{In the following we regard any
real function $x(\cdot)$ defined on a subset $D$ of $[0,\infty)$ as an element
of $\mathbb{R}^{[0,\infty)}$ by setting $x(t)=0$ for $t\notin D$.} \vskip0.2cm
Next, we introduce two useful dual operators. \newline

\begin{itemize}
\item For $T\in(0,\infty)$ let $G(t)=G(t,\cdot):\mathcal{\bar{S}}%
[0,T]\mapsto\mathbb{R}$ be a bounded linear operator on $\mathcal{\bar{S}%
}[0,T]$ for each $t$, uniformly bounded in $t\in\lbrack0,T]$. Then the map
\[
Y\mapsto\mathbb{E}[%
{\textstyle\int_{0}^{T}}
\left\langle G(t),Y_{t}\right\rangle dt];\quad Y\in\mathcal{\bar{S}}[0,T]
\]
is a bounded linear functional on the Hilbert space $\mathcal{\bar{S}}[0,T]$.
Therefore, by the Riesz representation theorem there exists a unique process
denoted by $G^{\ast}\in\mathcal{\bar{S}}[0,T],$ such that
\begin{equation}
{\mathbb{E}}[{%
{\textstyle\int_{0}^{T}}
}\left\langle G(t),Y_{t}\right\rangle dt]={\mathbb{E}}[{%
{\textstyle\int_{0}^{T}}
}G^{\ast}(t)Y(t)dt], \label{eq6.7a}%
\end{equation}
for all $Y\in\mathcal{\bar{S}}[0,T]$.


\item Proceeding as above, we also see that if $G_{\bar{m}}(t,\cdot
):[0,T]\times\mathcal{M}_{0,t}\mapsto L^{1}(\mathbb{P})$ is a bounded linear
operator on $\mathcal{M}_{0,t}$ for each $t$, uniformly bounded in $t$, then
the map
\[
M(\cdot)\mapsto%
{\textstyle\int_{0}^{T}}
\left\langle G_{\bar{m}}(t),M_{t}\right\rangle dt;\quad M_{t}=\mathcal{L}%
(X_{t})
\]
is a bounded linear functional on $\mathcal{M}_{0,t}$. Therefore, there exists
a unique process denoted by $G_{\bar{m}}^{\ast}(t)\in\mathcal{M}_{0,t}$ such
that
\begin{equation}
{%
{\textstyle\int_{0}^{T}}
}\left\langle G_{\bar{m}}(t),M_{t}\right\rangle dt={%
{\textstyle\int_{0}^{T}}
}G_{\bar{m}}^{\ast}(t)M(t)dt, \label{eq6.9a}%
\end{equation}
for all $M\in\mathcal{M}_{0,t}$.
\end{itemize}

We illustrate these operators by some auxiliary results.

\begin{lemma}
\label{lem} Consider the case when $G(t,\cdot):\mathcal{\bar{S}}%
[0,T]\mapsto\mathcal{\bar{S}}[0,T]$ has the form
\[
G(t,\bar{x})=\left\langle F,\bar{x}\right\rangle p(t)\text{, with }p\in
L_{0}^{2}.
\]
Then
\begin{equation}
G^{\ast}(t):=\left\langle F,p^{t}\right\rangle \label{eq4.8}%
\end{equation}
satisfies \eqref{eq6.7a}, where $p^{t}:=\{p(t+r)\}_{r\in\lbrack0,t]}$.
\end{lemma}

\noindent{Proof.} \quad\quad We must verify that if we define $G^{*}(t)$ by
\eqref{eq4.8}, then \eqref{eq6.7a} holds. To this end, choose $Y\in
\mathcal{S}_{\bar{x}}$ and consider%

\begin{align}
&
{\textstyle\int_{0}^{T}}
\left\langle F,p^{t}\right\rangle Y(t)dt={%
{\textstyle\int_{0}^{T}}
}\left\langle F,\{p(t+r)\}_{r\in\lbrack0,t]}\right\rangle Y(t)dt\nonumber\\
&  =%
{\textstyle\int_{0}^{T}}
\left\langle F,\{Y(t)p(t+r)\}_{r\in\lbrack0,t]}\right\rangle dt=\left\langle
F,\left\{
{\textstyle\int_{r}^{T+r}}
Y(u-r)p(u)du\right\}  _{r\in\lbrack0,t]}\right\rangle \nonumber\\
&  =\left\langle F,\left\{
{\textstyle\int_{0}^{T}}
Y(u-r)p(u)du\right\}  _{r\in\lbrack0,t]}\right\rangle ={%
{\textstyle\int_{0}^{T}}
}\left\langle F,Y_{u}\right\rangle p(u)du\nonumber\\
&  ={%
{\textstyle\int_{0}^{T}}
}\left\langle G(u),Y_{u}\right\rangle du.\nonumber
\end{align}
\hfill$\square$

\begin{example}
\label{exp} (i) For example, if $a\in\mathbb{R}^{[0,\infty)}$ is a bounded
function with compact support, let $F(\bar{x})$ be the averaging operator
defined by
\[
F(\bar{x})=\left\langle F,\bar{x}\right\rangle ={%
{\textstyle\int_{0}^{\infty}}
}a(r)x(r)dr
\]
when $\bar{x}=\{x(s)\}_{s\in\lbrack0,\infty)}$, then
\[
\left\langle F,p^{t}\right\rangle ={%
{\textstyle\int_{0}^{\infty}}
}a(r)p(t+r)dr.
\]
(ii) Similarly, if $F$ is evaluation at $t_{0}$, i.e.
\[
F(\bar{x})=x(t_{0})\text{ when }\bar{x}=\{x(s)\}_{s\in\lbrack0,\infty)},
\]
then
\[
\left\langle F,p^{t}\right\rangle =p(t+t_{0}).
\]

\end{example}

\color{black}

\section{The finite horizon case}

\noindent In this section we consider the case with a finite time horizon
$T<\infty$.\newline We are interested in the mean-field stochastic control
problem with elephant memory, composed of a controlled diffusion equation
defining the dynamics which are defined through the following equation:
\begin{equation}
\label{eq3.1a}\left\{
\begin{array}
[c]{ll}%
dX^{u}(t) & =b(t,X^{u}(t),X_{t}^{u},M^{u}(t),M_{t}^{u},u(t))dt+\sigma
(t,X^{u}(t),X_{t}^{u},M^{u}(t),M_{t}^{u},u(t))dB(t)\\
& +{%
{\textstyle\int_{\mathbb{R}_{0}}}
}\gamma(t,X^{u}(t),X_{t}^{u},M^{u}(t),M_{t}^{u},u(t),\zeta)\tilde{N}%
(dt,d\zeta);t\in\lbrack0,T],\\
X^{u}(0) & =x_{0},
\end{array}
\right.
\end{equation}
where $x_{0}\in%
\mathbb{R}
^{d}$ is a constant and $u\in\mathcal{U}_{adm}$ (the set of admissible
controls) is our control process, and with coefficients $b:[0,T]\times%
\mathbb{R}
^{d}\times C^{d}\times\mathcal{M}_{0}\times\mathcal{M}_{0,t}\times
U\rightarrow%
\mathbb{R}
^{d}$, $\sigma:[0,T]\times%
\mathbb{R}
^{d}\times C^{d}\times\mathcal{M}_{0}\times\mathcal{M}_{0,t}\times
U\rightarrow%
\mathbb{R}
^{d\times n}$ and $\gamma:[0,T]\times%
\mathbb{R}
^{d}\times C^{d}\times\mathcal{M}_{0}\times\mathcal{M}_{0,t}\times U\times%
\mathbb{R}
_{0}^{k}\rightarrow%
\mathbb{R}
^{k\times n}$ satisfying suitable assumptions (see below). Here and in the
following $U$ is the set of possible control values. For given $u\in
\mathcal{U}_{adm}$ we define its corresponding performance functional $J(u)$
by
\begin{equation}
J(u)=\mathbb{E}[{%
{\textstyle\int_{0}^{T}}
}f(t,X^{u}(t),X_{t}^{u},M^{u}(t),M_{t}^{u},u(t))dt+g(X^{u}(T),M^{u}(T))],
\label{eq3.4}%
\end{equation}
where $f:[0,T]\times%
\mathbb{R}
^{d}\times C^{d}\times\mathcal{M}_{0}\times\mathcal{M}_{0,t}\times
U\rightarrow%
\mathbb{R}
^{d}$ and $g:%
\mathbb{R}
^{d}\times\mathcal{M}_{0}\rightarrow%
\mathbb{R}
^{d}.$

\noindent We assume that $f(t,x,\bar{x},m,\bar{m},u)$ and $g(x,m)$ are
$\mathcal{F}_{t}$- and $\mathcal{F}_{T}$- measurable, respectively.
\vskip 0.2cm \noindent We consider the following \emph{finite horizon
mean-field elephant memory control problem}:

\begin{problem}
\label{prob-finite} Find $\hat{u}\in\mathcal{U}_{adm}$ such that
\[
J(\hat{u})=\sup_{u\in\mathcal{U}_{adm}}J(u).
\]

\end{problem}

\noindent For simplicity (but without loss of generality), from now on we will
consider only the one-dimensional case.\newline

\subsection{Existence and uniqueness of the MF-SDE with elephant memory}

\noindent We begin with the existence and uniqueness results for MF-SDE with
elephant memory. Consider the following equation for $X(t)=X^{\bar{u}}(t)$,
for fixed $\bar{u}\in\mathcal{U}_{adm}:$
\begin{equation}
\left\{
\begin{array}
[c]{ll}%
dX(t) & =b(t,X(t),X_{t},M(t),M_{t})dt+\sigma(t,X(t),X_{t},M(t),M_{t})dB(t)\\
& +%
{\textstyle\int_{\mathbb{R}_{0}}}
\gamma(t,X(t),X_{t},M(t),M_{t},\zeta)\tilde{N}(dt,d\zeta);t\in\lbrack0,T],\\
X(0) & =x_{0}.
\end{array}
\right.  \label{sde_fini}%
\end{equation}
We make the following assumptions on the coefficients $b:[0,T]\times%
\mathbb{R}
\times C\times\mathcal{M}_{0}\times\mathcal{M}_{0,t}\rightarrow%
\mathbb{R}
$, $\sigma:[0,T]\times%
\mathbb{R}
\times C\times\mathcal{M}_{0}\times\mathcal{M}_{0,t}\rightarrow%
\mathbb{R}
$ and $\gamma:[0,T]\times%
\mathbb{R}
\times C\times\mathcal{M}_{0}\times\mathcal{M}_{0,t}\times%
\mathbb{R}
_{0}\rightarrow%
\mathbb{R}
:$

\noindent Here the drift $b$, the volatility $\sigma$ and the jump coefficient
$\gamma$ are supposed to be $\mathbb{F}$-predictable.

\begin{description}
\item[(i)] The coefficients $b,$ $\sigma$ and $\gamma$ are Borel measurable.

\item[(ii)] There is a constant $C_{0}$ such that, for all $t\in
\lbrack0,T],\,\psi,\psi^{\prime}\in%
\mathbb{R}
,$ $\bar{\psi},\bar{\psi}^{\prime}\in C$, $m,m^{\prime}\in\mathcal{M}_{0}$ and
all $\bar{m},\bar{m}^{\prime}\in\mathcal{M}_{0,t}$, the following holds for
$h=b$ and for $h=\sigma$:
\end{description}%

\[
\left\{
\begin{array}
[c]{l}%
h\text{ is adapted and } \left\vert h(t,\psi,\bar{\psi},m,\bar{m})\right\vert
\leq C_{0},\\
h(t,\cdot,\cdot,\cdot,\cdot)\text{ is Lipschitz uniformly with respect to
}t,\\
\left\vert h\left(  t,\psi,\bar{\psi},m,\bar{m}\right)  -h\left(
t,\psi^{\prime},\bar{\psi}^{\prime},m^{\prime},\bar{m}^{\prime}\right)
\right\vert ^{2}\leq C_{0}(|\psi-\psi^{\prime}|^{2}+\underset{0\leq s\leq
t}{\sup}\left\vert \psi(s)-\psi^{\prime}(s)\right\vert ^{2}\\
\text{
\ \ \ \ \ \ \ \ \ \ \ \ \ \ \ \ \ \ \ \ \ \ \ \ \ \ \ \ \ \ \ \ \ \ \ \ }%
+||M(t)-M(t)||_{\mathcal{M}_{0}}^{2}+||M_{t}-M_{t}^{\prime}||_{\mathcal{M}%
_{0},t}^{2}).
\end{array}
\right.
\]
Similarly, we assume that $\gamma$ is predictable and
\[%
{\textstyle\int_{\mathbb{R}_{0}}}
\left\vert \gamma\left(  t,\psi,\bar{\psi},m,\bar{m},\zeta\right)  \right\vert
\nu(d\zeta)\leq C_{0},
\]%
\[%
\begin{array}
[c]{l}%
{\textstyle\int_{\mathbb{R}_{0}}}
\left\vert \gamma\left(  t,\psi,\bar{\psi},m,\bar{m},\zeta\right)
-\gamma\left(  t,\psi^{\prime},\bar{\psi}^{\prime},m^{\prime},\bar{m}^{\prime
},\zeta\right)  \right\vert ^{2}\nu(d\zeta)\leq C_{0}(|\psi-\psi^{\prime
}|+\underset{0\leq s\leq t}{\sup}\left\vert \psi(s)-\psi^{\prime
}(s)\right\vert ^{2}\\
\text{
\ \ \ \ \ \ \ \ \ \ \ \ \ \ \ \ \ \ \ \ \ \ \ \ \ \ \ \ \ \ \ \ \ \ \ \ \ \ \ \ \ \ \ \ }%
+||M(t)-M(t)||_{\mathcal{M}_{0}}^{2}+||M_{t}-M_{t}^{\prime}||_{\mathcal{M}%
_{0,t}}^{2}).
\end{array}
\]

\begin{theorem}
\label{thm-ex} Under the assumptions $(i)-(ii)$ our elephant memory MF-SDE
\[
\left\{
\begin{array}
[c]{ll}%
dX(t) & =b(t,X(t),X_{t},M(t),M_{t})dt+\sigma(t,X(t),X_{t},M(t),M_{t})dB(t)\\
& \text{ \ \ \ \ \ }+%
{\textstyle\int_{\mathbb{R}_{0}}}
\gamma(t,X(t),X_{t},M(t),M_{t},\zeta)\tilde{N}(dt,d\zeta);t\in\lbrack0,T],\\
X(0) & =x_{0}.
\end{array}
\right.
\]
for any initial condition $x_{0}\in\mathbb{R}$ admits a unique solution $X\in
S^{p}[0,T].$
\end{theorem}

\noindent We recall the following inequality which will be useful for our proof.

\begin{lemma}
[Kunita's inequality \cite{K}]Suppose {\normalsize $p\geq2$ }and{\normalsize
\[
X(t)=x_{0}+%
{\textstyle\int_{0}^{t}}
b(s)ds+%
{\textstyle\int_{0}^{t}}
\sigma(s)dB(s)+%
{\textstyle\int_{0}^{t}}
{\textstyle\int_{\mathbb{R}_{0}}}
\gamma(s,\zeta)\tilde{N}(ds,d\zeta).
\]
}Then there exists a positive constant {\normalsize $C_{p,T},$ }(depending
only on {\normalsize $p,T$}) such that the following inequality
holds{\normalsize
\begin{align*}
\mathbb{E}[\underset{0\leq t\leq T}{\sup}|X(t)|^{p}]  &  \leq C_{p,T}%
(|x_{0}|^{p}+\mathbb{E[}%
{\textstyle\int_{0}^{t}}
\{\left\vert b(s)\right\vert ^{p}+\left\vert \sigma(s)\right\vert ^{p}\\
&  +%
{\textstyle\int_{\mathbb{R}_{0}}}
\left\vert \gamma(s,\zeta)\right\vert ^{p}\nu(d\zeta)+(%
{\textstyle\int_{\mathbb{R}_{0}}}
\left\vert \gamma(s,\zeta)\right\vert ^{2}\nu(d\zeta))^{\tfrac{p}{2}}\}ds]).
\end{align*}
}
\end{lemma}

\noindent{Proof of Theorem \ref{thm-ex}.}\newline\quad\textbf{Existence}. For
the convenience of the reader, but without loss of generality of the method,
we assume that $b=\sigma=0$, but we can get the same result by using Kunita's
inequality above for $b\neq0$ and $\sigma\neq0$. In the following we denote by
$C_{p}$ the constant that may change from line to line.\newline Choose
arbitrary $X^{0}(t)$ with corresponding $X_{t}^{0},M^{0}(t),M_{t}^{0}$ and
consider inductively the equation
\[
\left\{
\begin{array}
[c]{ll}%
X(0) & :=x_{0},\\
X^{n+1}(t) & =x_{0}+%
{\textstyle\int_{0}^{t}}
{\textstyle\int_{\mathbb{R}_{0}}}
\gamma(s,X^{n}(s),X_{s}^{n},M^{n}(s),M_{s}^{n},\zeta)\tilde{N}(ds,d\zeta
),\text{ }t\in\left[  0,T\right]  ,\,n\geq0.
\end{array}
\right.
\]
It is clear that $X^{n}(t)\in S^{p}\left[  0,T\right]  $, for all $n\geq0$.
Let $\overline{X}^{n}:=X^{n+1}-X^{n}.$ Then, by the Kunita's inequality ( for
$b=\sigma=0)$, the following estimation holds for all $p\geq2:$%
\[%
\begin{array}
[c]{l}%
\mathbb{E}[\left\vert \overline{X}^{n}(s)\right\vert ^{p}]\leq C_{p}%
(\mathbb{E}[%
{\textstyle\int_{0}^{t}}
{\textstyle\int_{\mathbb{R}_{0}}}
|\gamma(s,X^{n}(s),X_{s}^{n},M^{n}(s),M_{s}^{n},\zeta)\\
-\gamma(s,X^{n-1}(s),X_{s}^{n-1},M^{n-1}(s),M_{s}^{n-1},\zeta)|^{p}\nu
(d\zeta)ds]\\
+\mathbb{E}[%
{\textstyle\int_{0}^{t}}
(%
{\textstyle\int_{\mathbb{R}_{0}}^{2}}
|\gamma(s,X^{n}(s),X_{s}^{n},M^{n}(s),M_{s}^{n},\zeta)+\\
-\gamma(s,X^{n-1}(s),X_{s}^{n-1},M^{n-1}(s),M_{s}^{n-1},\zeta)|\nu
(d\zeta))^{\tfrac{p}{2}}ds]),\text{ }t\in\left[  0,T\right]  ,n\geq1.
\end{array}
\]
Applying the Lipschitz assumption (ii), we get
\[%
\begin{array}
[c]{l}%
\mathbb{E}[\underset{s\leq t}{\sup}\left\vert \overline{X}^{n}(s)\right\vert
^{p}]\leq C_{p}\mathbb{E}[\underset{s\leq t}{\sup}|\overline{X}^{n-1}%
(s)|+||M^{n}(s)-M^{n-1}(s)||_{\mathcal{M}_{0}}\\
+%
{\textstyle\int_{0}^{t}}
(\underset{r\leq s}{\sup}|\overline{X}^{n-1}(r)|+||M^{n}(r)-M^{n-1}%
(r)||_{\mathcal{M}_{0}})^{2}ds]^{p/2}\\
\leq C_{p}\mathbb{E[}\underset{s\leq t}{\sup}|\overline{X}^{n-1}%
(s)|^{p}],\text{ }t\in\left[  0,T\right]  ,\,n\geq1.
\end{array}
\]
Hence, from a standard argument we see that there is some $X\in\underset
{p>1}{\cap}S^{p}[0,T],$ such that
\[
\mathbb{E}[\underset{t\in\left[  0,T\right]  }{\sup}\left\vert X^{n}%
(t)-X(t)\right\vert ^{p}]\underset{n\rightarrow\infty}{\rightarrow}0,\text{
for all }p\geq2.
\]
Finally, taking the limit in the Picard iteration as $n\rightarrow+\infty,$
yields
\[
X(t)=x_{0}+%
{\textstyle\int_{0}^{t}}
{\textstyle\int_{\mathbb{R}_{0}}}
\gamma(s,X(s),X_{s},M(s),M_{s},u_{s},\zeta)\tilde{N}(ds,d\zeta),\text{ }%
t\in\left[  0,T\right]  .
\]
\textbf{Uniqueness.} The proof of uniqueness is obtained by the estimate of
the difference of two solutions, and it is carried out similarly to the
argument above. \hfill$\square$ \bigskip

\subsection{Stochastic maximum principles}

We now turn to the problem of optimal control of the mean-field equation
\eqref{eq3.1a} with performance functional \eqref{eq3.4}. Because of the
mean-field terms, it is natural to consider the two-dimensional system
$(X(t),M(t))$, where the dynamics for $M(t)$ is the following:
\[%
\begin{cases}
dM(t)=\beta(M(t)dt,\\
M(0)\in\mathcal{M}_{0},
\end{cases}
\]
where we have put $\beta(M(t))=M^{\prime}(t)$. See Lemma \ref{2.6}.

\noindent Let $\mathcal{R}$ denote the set of Borel measurable functions $r:%
\mathbb{R}
_{0}\rightarrow%
\mathbb{R}
.$\newline Define the Hamiltonian $H:[0,T]\times%
\mathbb{R}
\times C\times\mathcal{M}_{0}\times\mathcal{M}_{0,t}\times U\times%
\mathbb{R}
\times%
\mathbb{R}
\times\mathcal{R}\times C_{a}([0,T],\mathcal{M}_{0})\rightarrow%
\mathbb{R}
$, as follows%

\begin{align*}
H(t,x,\bar{x},m,\bar{m},u,p^{0},q^{0},r^{0},p^{1})  &  :=f(t,x,\bar{x}%
,m,\bar{m},u)+p^{0}b(t,x,\bar{x},m,\bar{m},u)\\
+q^{0}\sigma(t,x,\bar{x},m,\bar{m},u)  &  +%
{\textstyle\int_{\mathbb{R}_{0}}}
r^{0}(\zeta)\gamma(t,x,\bar{x},m,\bar{m},u,\zeta)\nu(d\zeta)+\left\langle
p^{1},\beta(m)\right\rangle ;\quad t\in\lbrack0,T],
\end{align*}
and $H(t,x,\bar{x},m,\bar{m},u,p^{0},q^{0},r^{0},p^{1})=0$ for all $t>T.$

\noindent We assume that all the coefficients $f,b,\sigma,\gamma$ and $g$ are
continuously differentiable $(C^{1})$ with respect to $x,$ $u,$ and admit
Fréchet derivatives with respect to $\overline{x},m,\overline{m}$. Then the
same holds for the Hamiltonian $H$.\newline

\noindent We define the adjoint processes $(p^{0},q^{0},r^{0}),(p^{1}%
,q^{1},r^{1})$ as the solution of the following finite horizon backward
stochastic differential equations (BSDEs):%

\begin{equation}
\left\{
\begin{array}
[c]{ll}%
dp^{0}(t) & =-\Big\{\tfrac{\partial H}{\partial x}(t)+\mathbb{E[}\nabla
_{\bar{x}}^{\ast}H(t)|\mathcal{F}_{t}]\Big\}dt+q^{0}(t)dB(t)\\
& \quad+{{{%
{\textstyle\int_{\mathbb{R}_{0}}}
}}}r^{0}(t,\zeta)\tilde{N}(dt,d\zeta);\quad t\in\lbrack0,T],\\
p^{0}(t) & =\tfrac{\partial g}{\partial x}(T);\quad t\geq T,\\
q^{0}(t) & =r^{0}(t,\cdot)=0;\quad t>T,
\end{array}
\right.  \label{bsde0}%
\end{equation}%
\begin{equation}
\left\{
\begin{array}
[c]{ll}%
dp^{1}(t) & =-\left(  \nabla_{m}H(t)+\mathbb{E}\left[  \nabla_{\bar{m}}^{\ast
}H(t)|\mathcal{F}_{t}\right]  \right)  dt+q^{1}(t)dB(t)\\
& \quad+\int_{\mathbb{R}_{0}}r^{1}(t,\zeta)\tilde{N}(dt,d\zeta);\quad
t\in\lbrack0,T],\\
p^{1}(t) & =\nabla_{m}g(T);\quad t\geq T,\\
q^{1}(t) & =r^{1}(t,\cdot)=0;\quad t>T,
\end{array}
\right.  \label{eqp1}%
\end{equation}
where $g(T)=g(X(T),M(T))$ and
\[
H(t)=H(t,x,\bar{x},m,\bar{m},u,p^{0},q^{0},r^{0},p^{1})_{x=X(t),\bar{x}%
=X_{t},m=M(t),\bar{m}=M_{t},u=u(t),p^{0}=p^{0}(t),q^{0}=q^{0}(t),r^{0}%
=r^{0}(t,\zeta),p^{1}=p^{1}(t)}.
\]

\noindent In the next section, we will give an example on how to calculate
this adjoint operator in particular cases.

\noindent We are now able to give a sufficient (a verification theorem) and a
necessary maximum principle. \vskip 0.2cm  \emph{We do not give
the proof of the following result, since it is similar to the proof in the
infinite horizon case, which will be discussed in Section 4.} \color{black}

\begin{theorem}
[Sufficient conditions of optimality]Let $\hat{u}\in\mathcal{U}_{adm}$ with
corresponding solutions $\hat{X}$ and $(\hat{p}^{0},\hat{q}^{0},\hat{r}^{0})$
and $(\hat{p}^{1},\hat{q}^{1},\hat{r}^{1})$ of the forward and the backward
stochastic differential equations \eqref{sde_fini}, \eqref{bsde0}
and (\ref{eqp1}) respectively. \color{black} Suppose that

\begin{enumerate}
\item (Concavity) The Hamiltonian is such that
\[
(x,\overline{x},m,\bar{m},u)\mapsto H(t,x,\bar{x},m,\bar{m},u,\hat{p}%
^{0}(t),\hat{q}^{0}(t),\hat{r}^{0}(t,\zeta),\hat{p}^{1}(t),\omega),
\]
and the terminal condition
\[
(x,m)\mapsto g(x,m,\omega),
\]
are concave $\mathbb{P}$-a.s. for each $t$.

\item (Maximum condition)
\begin{align}
&  \mathbb{E}[H(t,\hat{X}(t),\hat{X}_{t},\hat{M}(t),\hat{M}_{t},\hat
{u}(t),\hat{p}^{0}(t),\hat{q}^{0}(t),\hat{r}^{0}(t,\cdot),\hat{p}%
^{1}(t))|\mathcal{G}_{t}]\nonumber\\
&  =\sup_{u\in\mathcal{U}_{adm}}\mathbb{E}[H(t,\hat{X}(t),\hat{X}_{t},\hat
{M}(t),\hat{M}_{t},u,\hat{p}^{0}(t),\hat{q}^{0}(t),\hat{r}^{0}(t,\cdot
),\hat{p}^{1}(t))|\mathcal{G}_{t}],
\end{align}
$\mathbb{P}$-a.s. for each $t\in\lbrack0,T].$
\end{enumerate}

Then $\hat{u}$ is an optimal control for Problem \ref{prob-finite}.
\end{theorem}

\vskip 0.3cm \noindent Next we consider a converse, in the sense that we look
for necessary conditions of optimality. To this end, we make the following assumptions:

\begin{itemize}
\item \emph{Assumption A1}.\newline Whenever $u\in\mathcal{U}_{adm}$, and
$\pi\in\mathcal{U}_{adm}$ is bounded, there exists $\epsilon>0$ such that for
$\lambda\in(-\epsilon,\epsilon)$ we have
\[
u+\lambda\pi\in\mathcal{U}_{adm}.
\]

\item \emph{Assumption A2}.\newline For each $t_{0}\in\left[  0,T\right]  $
and each bounded $\mathcal{G}_{t_{0}}$-measurable random variables $\alpha,$
the process
\[
\pi(t)=\alpha\mathbf{1}_{\left(  t_{0},T\right]  }(t)
\]
belongs to $\mathcal{U}_{adm}$.

\item \emph{Assumption A3}.\newline In general, if $K^{u}=(K^{u}%
(t))_{t\in\lbrack0,T]}$ is a process depending on $u$, and if $\pi
\in\mathcal{U}$ we define the operator $D=D_{\pi}$ on $K$ by
\[
DK^{u}(t):=D_{\pi}K^{u}(t)=\tfrac{d}{d\lambda}K^{u+\lambda\pi}(t)|_{\lambda
=0},
\]
whenever the derivative exists. In particular, we define the \emph{derivative
process} $Z=Z_{\pi}=(Z(t))_{t\in\lbrack0,T]}$ by
\[
Z(t)=DX^{u}(t):=\tfrac{d}{d\lambda}X^{u+\lambda\pi}(t)|_{\lambda=0}.
\]

We assume that for all bounded $\pi\in\mathcal{U}_{adm}$ the derivative
process $Z(t)=Z_{\pi}(t)$ exists and satisfies the equation
\begin{equation}
\left\{
\begin{array}
[c]{ll}%
dZ(t) & =\left[  \tfrac{\partial b}{\partial x}(t)Z(t)+\langle\nabla
_{\overline{x}}b(t),Z_{t}\rangle+\langle\nabla_{m}b(t),DM(t)\rangle\right. \\
& +\left.  \langle\nabla_{\overline{m}}b(t),DM_{t}\rangle+\frac{\partial
b}{\partial u}(t)\pi(t)\right]  dt\\
& +\left[  \tfrac{\partial\sigma}{\partial x}(t)Z(t)+\langle\nabla
_{\overline{x}}\sigma(t),Z_{t}\rangle+\langle\nabla_{m}\sigma(t),DM(t)\rangle
\right. \\
& \left.  +\langle\nabla_{\overline{m}}\sigma(t),DM_{t}\rangle+\frac
{\partial\sigma}{\partial u}(t)\pi(t)\right]  dB(t)\\
& +%
{\textstyle\int_{\mathbb{R}_{0}}}
\left[  \tfrac{\partial\gamma}{\partial x}(t,\zeta)Z(t)+\langle\nabla
_{\overline{x}}\gamma(t,\zeta),Z_{t}\rangle+\langle\nabla_{m}\gamma
(t,\zeta),DM(t)\rangle\right. \\
& +\left.  \langle\nabla_{\overline{m}}\gamma(t,\zeta),DM_{t}\rangle
+\frac{\partial\gamma}{\partial u}(t,\zeta)\pi(t)\right]  \tilde{N}%
(dt,d\zeta);\quad t\in\lbrack0,T],\\
Z(0) & =0.
\end{array}
\right.  \label{eq3.25}%
\end{equation}

\end{itemize}

\begin{remark}
Using the Itô formula we see that Assumption A3 holds under reasonable
smoothness conditions on the coefficients of the equation. A proof for a
similar system is given in Lemma 12 in Agram and Øksendal \cite{AO3}. We omit
the details.
\end{remark}

\emph{We do not give the proof of the following result, since it is similar to
the proof in the infinite horizon case, which will be discussed in Section 4.} 

\begin{theorem}
\label{nece-fin} Let $\hat{u}\in\mathcal{U}_{adm}$ with corresponding
solutions $\hat{X}$ and $(\hat{p}^{0},\hat{q}^{0},\hat{r}^{0})$ and $(\hat
{p}^{1},\hat{q}^{1},\hat{r}^{1})$ of the forward and the backward stochastic
differential equations \eqref{sde_fini}, \eqref{bsde0} and
(\ref{eqp1}) \color{black} respectively with corresponding derivative process
$\hat{Z}$ given by \eqref{eq3.25}.

Then the following are equivalent:

\begin{itemize}
\item
\begin{equation}
\label{eq3.26}\tfrac{d}{d\lambda}J(\hat{u}+\lambda\pi)|_{\lambda=0}=0 \text{
for all bounded } \pi\in\mathcal{U}_{adm}.
\end{equation}

\item
\begin{equation}
\mathbb{E[}\tfrac{\partial H}{\partial u}(t,\hat{X}(t),\hat{X}_{t},\hat
{M}(t),\hat{M}_{t},u,\hat{p}(t),\hat{q}(t),\hat{r}(t,\cdot))_{u=\hat{u}%
}|\mathcal{G}_{t}]=0. \label{eq3.27}%
\end{equation}

\end{itemize}
\end{theorem}

\subsection{Example: A mean-field LQ control problem}

\noindent As an example, consider the following optimization problem which is
to maximize the performance functional
\[
J(u)=\mathbb{E}[-\tfrac{1}{2}X^{2}(T)-\tfrac{1}{2}%
{\textstyle\int_{0}^{T}}
u^{2}(t)dt],
\]
where $X(t)$ is subject to%
\begin{equation}
\left\{
\begin{array}
[c]{ll}%
dX(t) & =\mathbb{E}[X(t)](b_{0}+u(t))dt+\sigma_{0}\mathbb{E}[X(t)]dB(t)\\
& +%
{\textstyle\int_{\mathbb{R}_{0}}}
\gamma_{0}(\zeta)\mathbb{E}[X(t)]\tilde{N}(dt,d\zeta),\\
X(0) & =x_{0}\in\mathbb{R},
\end{array}
\right.  \label{eq4.20}%
\end{equation}
for some given constants $b_{0},$ $\sigma_{0}$ and $\gamma_{0}(\zeta)>-1$
a.s. $\nu$.

\noindent We associate to this problem the Hamiltonian%
\begin{align}
H(t,m,u,p^{0},q^{0},r^{0},p^{1})  &  =-\tfrac{1}{2}u^{2}+F(m)(b_{0}%
+u)p^{0}+F(m)\sigma_{0}q^{0}\label{eq4.21}\\
&  +{%
{\textstyle\int_{\mathbb{R}_{0}}}
}F(m)\gamma_{0}(\zeta)r^{0}(\zeta)\nu(d\zeta)+\left\langle p^{1}%
,\beta(m)\right\rangle .\nonumber
\end{align}
Here
\begin{align*}
&  b(t,X(t),X_{t},M(t),M_{t})=F(M(t))(b_{0}+u(t)),\\
&  \sigma(t,X(t),X_{t},M(t),M_{t})=F(M(t))\sigma_{0},\\
&  \gamma(t,X(t),X_{t},M(t),M_{t},\zeta)=F(M(t))\gamma_{0}(\zeta)r(\zeta
)\nu(d\zeta),
\end{align*}
where the operator $F$ is defined by
\[
F(m)=%
{\textstyle\int_{\mathbb{R}}}
xdm(x);\quad m\in\mathcal{M}_{0},
\]
so that
\[
F(M(t))=%
{\textstyle\int_{\mathbb{R}}}
xdM(t)(x)=\mathbb{E}[X(t)]\text{ when }M(t)=\mathcal{L}(X(t)).
\]

\noindent Note that, since $H$ does not depend on $x$, $\bar{x}$, $\bar{m}$,
we have
\[
\tfrac{\partial H}{\partial x}(t)=\nabla_{\bar{x}}H(t)=\nabla_{\bar{m}%
}H(t)=0.
\]
And, since $m\mapsto F(m)$ and $m\mapsto\beta(m)$ are linear, we have
\[
\nabla_{m}H(t)=F(\cdot)(b_{0}+u)p^{0}(t)+F(\cdot)\sigma_{0}q^{0}%
(t)\label{eq4.21}\newline+{\textstyle\int_{\mathbb{R}_{0}}}F(\cdot)\gamma
_{0}(\zeta)r^{0}(\zeta)\nu(d\zeta)+\left\langle p^{1},\beta(\cdot
)\right\rangle .
\]
Hence, the adjoint equation for $(p^{0},q^{0},r^{0})$ is
\begin{equation}
\left\{
\begin{array}
[c]{ll}%
dp^{0}(t) & =q^{0}(t)dB(t)+%
{\textstyle\int_{\mathbb{R}_{0}}}
r^{0}(t,\zeta)\tilde{N}(dt,d\zeta);\quad0\leq t\leq T,\\
p^{0}(T) & =-X(T),
\end{array}
\right.  \label{lbsde}%
\end{equation}
and the adjoint equation for $(p^{1},q^{1},r^{1})$ is%
\[
\left\{
\begin{array}
[c]{ll}%
dp^{1}(t) & =-[F(\cdot)(b_{0}+u)p^{0}(t)+F(\cdot)\sigma_{0}q^{0}%
(t)+F(\cdot)\gamma_{0}(\zeta)r^{0}(\zeta)\nu(d\zeta)+<p^{1},\beta(\cdot)>]dt\\
& +q^{1}(t)dB(t)+\int_{\mathbb{R}_{0}}r^{1}(t,\zeta)\tilde{N}(dt,d\zeta);\quad
t\in\lbrack0,T],\\
p^{1}(T) & =0.
\end{array}
\right.
\]

\noindent The map $u\mapsto H(u)$ is maximal when $\tfrac{\partial H}{\partial
u}=0$, i.e., when
\begin{equation}
u=\hat{u}(t)=\mathbb{E}[\hat{X}(t)]\hat{p}^{0}(t)=-\mathbb{E}[\hat
{X}(t)]\mathbb{E}[\hat{X}(T)|\mathcal{F}_{t}]. \label{eq3.14}%
\end{equation}
Substituting this into \eqref{eq4.20} we get that $Y(t):=\mathbb{E}[\hat
{X}(t)]$ satisfies the following Riccati equation
\begin{equation}
\left\{
\begin{array}
[c]{ll}%
Y^{\prime}(t) & =b_{0}Y(t)-Y^{2}(t)Y(T);\quad0\leq t\leq T,\\
Y(0) & =x_{0}.
\end{array}
\right.  \label{eq4.31}%
\end{equation}
Solving this Riccati equation, we find an explicit expression for $Y(t)$ in
terms of $Y(T)$ and hence by putting $t=T$ also an explicit expression for
$Y(T)$, and then we find $Y(t)$ for all $t\in\lbrack0,T].$

\noindent Equation (\ref{eq4.31}) has the solution:%

\[
Y(t)=\mathbb{E}[\hat{X}(t)]=\tfrac{b_{0}x_{0}\exp(b_{0}t)}{(b_{0}%
-x_{0}\mathbb{E}[\hat{X}(T)])(1+\exp(b_{0}t))}.
\]
Consequently,%
\[
Y(T)=\mathbb{E}[\hat{X}(T)]=\tfrac{b_{0}x_{0}\exp(b_{0}T)}{(b_{0}%
-x_{0}\mathbb{E}[\hat{X}(T)])(1+\exp(b_{0}T))}.
\]

\noindent Then we see that we also know $\mathbb{E}[\hat{X}(T)|\mathcal{F}%
_{t}]$ by the equation
\begin{align*}
\mathbf{K}(t)  &  =\mathbb{E}[\hat{X}(T)|\mathcal{F}_{t}]\\
&  =\mathbf{K}(0)+%
{\textstyle\int_{0}^{t}}
Y(s)\sigma_{0}dB(s)+%
{\textstyle\int_{0}^{t}}
{\textstyle\int_{\mathbb{R}_{0}}}
Y(s)\gamma_{0}(\zeta)\tilde{N}(ds,d\zeta).
\end{align*}
Thus we have proved the following:

\begin{theorem}
The optimal control $\hat{u}$ of the mean-field LQ problem is given by
\[
\hat{u}(t)=-\mathbb{E}[\hat{X}(t)]\mathbb{E}[\hat{X}(T)|\mathcal{F}_{t}],
\]
with $\mathbb{E}[\hat{X}(t)]$ and $\mathbb{E}[\hat{X}(T)|\mathcal{F}_{t}]$
given above.
\end{theorem}

\section{The infinite horizon case}

\noindent We now study the case when the time horizon is $[0,\infty)$.
Consider the equation
\begin{equation}
\left\{
\begin{array}
[c]{ll}%
dX(t) & =b(t,X(t),X_{t},M(t),M_{t},u(t))dt+\sigma(t,X(t),X_{t},M(t),M_{t}%
,u(t))dB(t)\\
& +%
{\textstyle\int_{\mathbb{R}_{0}}}
\gamma(t,X(t),X_{t},M(t),M_{t},u(t),\zeta)\tilde{N}(dt,d\zeta);t\in
\lbrack0,\infty),\\
X(0) & =x_{0},
\end{array}
\right.  \label{F}%
\end{equation}
where $x_{0}\in%
\mathbb{R}
$ is the initial condition, $u\in\mathcal{U}_{adm}$, and the coefficients
$b:[0,\infty)\times%
\mathbb{R}
\times C\times\mathcal{M}_{0}\times\mathcal{M}_{0,t}\times U\rightarrow%
\mathbb{R}
$, $\sigma:[0,\infty)\times%
\mathbb{R}
\times C\times\mathcal{M}_{0}\times\mathcal{M}_{0,t}\times U\rightarrow%
\mathbb{R}
$ and $\gamma:[0,\infty)\times%
\mathbb{R}
\times C\times\mathcal{M}_{0}\times\mathcal{M}_{0,t}\times U\times%
\mathbb{R}
_{0}\rightarrow%
\mathbb{R}
$ are $\mathcal{F}_{t}$-measurable. Here $C$ stands for the space of
$\mathbb{R}$-valued continuous functions defined over the time interval
$[0,\infty).$ We assume that%
\[
\mathbb{E}[{%
{\textstyle\int_{0}^{\infty}}
|X(s)|}^{2}ds]<\infty.
\]
\newline For given $u\in\mathcal{U}_{adm},$ we define its corresponding
performance functional by
\begin{equation}
J(u)=\mathbb{E}[{%
{\textstyle\int_{0}^{\infty}}
}f(t,X(t),X_{t},M(t),M_{t},u(t))dt], \label{P}%
\end{equation}
where the reward function $f:[0,\infty)\times%
\mathbb{R}
\times C\times\mathcal{M}_{0}\times\mathcal{M}_{0,t}\times U\rightarrow%
\mathbb{R}
$ is assumed to satisfy the condition
\[
\mathbb{E}[{%
{\textstyle\int_{0}^{\infty}}
|}f(t,X(t),X_{t},M(t),M_{t},u(t))|^{2}dt]<\infty,\quad\text{ for all }%
u\in\mathcal{U}_{adm}.
\]
We consider the following \emph{infinite horizon mean-field elephant memory
control problem}:

\begin{problem}
\label{prob} Find $\hat{u}\in\mathcal{U}_{adm}$ such that
\[
J(\hat{u})=\sup_{u\in\mathcal{U}_{adm}}J(u).
\]

\end{problem}

\noindent Define the Hamiltonian $H:[0,\infty)\times%
\mathbb{R}
\times C\times\mathcal{M}_{0}\times\mathcal{M}_{0,t}\times U\times%
\mathbb{R}
\times%
\mathbb{R}
\times\mathcal{R}\times C_{a}([0,T],\mathcal{M}_{0})\rightarrow%
\mathbb{R}
$, by%

\begin{equation}%
\begin{array}
[c]{c}%
H(t,x,\bar{x},m,\bar{m},u,p^{0},q^{0},r^{0},p^{1}):=f(t,x,\bar{x},m,\bar
{m},u)+p^{0}b(t,x,\bar{x},m,\bar{m},u)\\
\text{ \ \ \ \ \ \ \ \ \ \ \ \ }+q^{0}\sigma(t,x,\bar{x},m,\bar{m},u)+%
{\textstyle\int_{\mathbb{R}_{0}}}
r^{0}(\zeta)\gamma(t,x,\bar{x},m,\bar{m},u,\zeta)\nu(d\zeta)+\left\langle
p^{1},m^{\prime}\right\rangle .
\end{array}
\label{eq3.6}%
\end{equation}
In the following we assume that all the coefficients $f,b,\sigma$ and $\gamma$
are continuously differentiable $(C^{1})$ with respect to $x$ and admit
Fréchet derivatives with respect to $\overline{x},m,\overline{m}$ and
$u$.\newline

\noindent Associated to the control $\hat{u}$ we define the following infinite
horizon BSDE for the adjoint processes $(\hat{p}^{0},\hat{q}^{0},\hat{r}^{0}%
)$, $(\hat{p}^{1},\hat{q}^{1},\hat{r}^{1})$:
\begin{equation}%
\begin{array}
[c]{ll}%
dp^{0}(t) & =-\Big\{\tfrac{\partial H}{\partial x}(t)+\mathbb{E[}\nabla
_{\bar{x}}^{\ast}H(t)|\mathcal{F}_{t}]\Big\}dt\\
& +q^{0}(t)dB(t)+{{{%
{\textstyle\int_{\mathbb{R}_{0}}}
}}}r^{0}(t,\zeta)\tilde{N}(dt,d\zeta);\quad t\geq0,
\end{array}
\label{eq4.11}%
\end{equation}%
\begin{equation}%
\begin{array}
[c]{cc}%
dp^{1}(t) & =-[\nabla_{m}H(t)+\mathbb{E[}\nabla_{\bar{m}}^{\ast}%
H(t)|\mathcal{F}_{t}]]dt+q^{1}(t)dB(t)\\
& +\int_{\mathbb{R}_{0}}r^{1}(t,\zeta)\tilde{N}(dt,d\zeta);\quad t\geq0.
\end{array}
\label{eq4.11'}%
\end{equation}

\begin{remark}
Note that without further conditions there are infinitely many solutions
$(\hat{p}^{0},\hat{q}^{0},\hat{r}^{0})$ and $(\hat{p}^{1},\hat{q}^{1},\hat
{r}^{1})$ of these equations.
\end{remark}

\subsection{Sufficient infinite horizon maximum principle}

\noindent In this subsection, we give sufficient conditions which ensure the
existence of an optimal control in the infinite horizon case.

\begin{theorem}
[Sufficient condition of optimality]Let $\hat{u}\in\mathcal{U}_{adm}$ with
corresponding solution $\hat{X}$ of the forward stochastic differential
equation \eqref{F}. Assume that $(\hat{p}^{0},\hat{q}^{0},\hat{r}^{0})$ and
$(\hat{p}^{1},\hat{q}^{1},\hat{r}^{1})$ is \emph{some} solution of the
associated backward stochastic differential equations \eqref{eq4.11} and
\eqref{eq4.11'} respectively. Suppose the following holds:

\begin{enumerate}
\item (Concavity) The function
\[
(x,\overline{x},m,\bar{m},u)\mapsto H(t,x,\bar{x},m,\bar{m},u,\hat{p}^{0}%
,\hat{q}^{0},\hat{r}^{0},\hat{p}^{1}),
\]
is concave $\mathbb{P}$-a.s. for each $t,\hat{p}^{0},\hat{q}^{0},\hat{r}%
^{0},\hat{p}^{1}$.

\item (Maximum condition)%
\begin{equation}
\mathbb{E[}\hat{H}(t)|\mathcal{G}_{t}]=\underset{u\in\mathcal{U}_{adm}}{\text{
}\sup}\mathbb{E}\left[  H(t)|\mathcal{G}_{t}\right]  , \label{maxQ}%
\end{equation}
$\mathbb{P}$-a.s. for each $t\geq0.$

\item (Transversality condition) For all $u\in\mathcal{U}_{adm}$ with
corresponding solution $X^{u}=X$ we have
\begin{equation}
\underset{T\rightarrow\infty}{\overline{\lim}}\mathbb{E}[\hat{p}%
^{0}(T)(X(T)-\hat{X}(T))]+\underset{T\rightarrow\infty}{\overline{\lim}%
}\mathbb{E}[\hat{p}^{1}(T)(M(T)-\hat{M}(T))]\geq0, \label{tcond_1}%
\end{equation}

\end{enumerate}

Then $\hat{u}$ is an optimal control for Problem \ref{prob}.
\end{theorem}

\noindent{Proof.} \quad Choose arbitrary $u\in\mathcal{U}_{adm}$. We want to
show that $J(u)\leq J(\hat{u}),$ i.e.,
\begin{equation}
A:=J(u)-J(\hat{u})=\mathbb{E[}{%
{\textstyle\int_{0}^{\infty}}
}\{f(t)-\hat{f}(t)\}dt]\leq0,\label{J2}%
\end{equation}
where we have used the simplified notation $\hat{f}(t):=f(t,\hat{X}(t),\hat
{X}_{t},\hat{M}(t),\hat{M}_{t},\hat{u}(t))$ and so on.\newline By concavity of
the Hamiltonian \eqref{eq3.6}, we have%
\begin{align}
A &  =\mathbb{E[}{%
{\textstyle\int_{0}^{\infty}}
}\{H(t)-\hat{H}(t)-\hat{p}^{0}(t)\tilde{b}(t)-\hat{q}^{0}(t)\tilde{\sigma}(t)-%
{\textstyle\int_{\mathbb{R}_{0}}}
\hat{r}^{0}(t,\zeta)\tilde{\gamma}(t,\zeta)\nu(d\zeta)\}dt]\label{eq3.13}\\
&  \leq\mathbb{E[}{%
{\textstyle\int_{0}^{\infty}}
}\{\tfrac{\partial\hat{H}}{\partial x}(t)\tilde{X}(t)+\langle\nabla_{\bar{x}%
}\hat{H}(t),\tilde{X}_{t}\rangle+\langle\nabla_{m}\hat{H}(t),\tilde
{M}(t)\rangle+\langle\nabla_{\bar{m}}\hat{H}(t),\tilde{M}_{t}\rangle
+\tfrac{\partial H}{\partial u}(t)\tilde{u}(t)\nonumber\\
&  -\hat{p}^{0}(t)\tilde{b}(t)-\hat{q}^{0}(t)\tilde{\sigma}(t)-%
{\textstyle\int_{\mathbb{R}_{0}}}
\hat{r}^{0}(t,\zeta)\tilde{\gamma}(t,\zeta)\nu(d\zeta)\}dt],\nonumber
\end{align}
where $\tilde{b}(t)=b(t)-\hat{b}(t)$, etc.\newline For fixed $T\geq0$, define
an increasing sequence of stopping times $\tau_{n}$, as follows%

\[%
\begin{array}
[c]{l}%
\tau_{n}(\cdot):=T\wedge\inf\{t\geq0:\int_{0}^{t}((\hat{p}^{0}(s)\tilde
{\sigma}(s))^{2}+(\hat{q}^{0}(s)\tilde{X}(s))^{2}\\
+%
{\textstyle\int_{\mathbb{R}_{0}}}
\{(\hat{r}^{0}(s,\zeta)\tilde{X}(s))^{2}+(\hat{p}^{0}(s)\tilde{\gamma}%
(s,\zeta))^{2}\}\nu(d\zeta))ds\geq n\},n\in%
\mathbb{N}
,
\end{array}
\]
it clearly holds that {\normalsize $\tau_{n}\rightarrow T$ }$\mathbb{P}$-a.s.
By the Itô formula applied to $\hat{p}^{0}({\normalsize \tau_{n}})\tilde
{X}({\normalsize \tau_{n}})$, we get
\begin{align*}
&  \mathbb{E}[\hat{p}^{0}(T)\tilde{X}(T)]=\underset{n\rightarrow\infty}{\lim
}\mathbb{E}[\hat{p}^{0}({\normalsize \tau_{n}})\tilde{X}({\normalsize \tau
_{n}})]\\
&  =\underset{n\rightarrow\infty}{\lim}\mathbb{E[}{%
{\textstyle\int_{0}^{{\normalsize \tau_{n}}}}
}\hat{p}^{0}(t)d\tilde{X}(t)+{%
{\textstyle\int_{0}^{{\normalsize \tau_{n}}}}
}\tilde{X}(t)d\hat{p}^{0}(t)+{%
{\textstyle\int_{0}^{{\normalsize \tau_{n}}}}
}\hat{q}^{0}(t)\tilde{\sigma}(t)dt\\
&  +{%
{\textstyle\int_{0}^{{\normalsize \tau_{n}}}}
}%
{\textstyle\int_{\mathbb{R}_{0}}}
\hat{r}(t,\zeta)\tilde{\gamma}(t,\zeta)\nu(d\zeta)dt]\\
&  =\underset{n\rightarrow\infty}{\lim}\mathbb{E[}{%
{\textstyle\int_{0}^{{\normalsize \tau_{n}}}}
}\{\hat{p}^{0}(t)\tilde{b}(t)-\tilde{X}(t)(\tfrac{\partial\hat{H}}{\partial
x}(t)+\nabla_{\bar{x}}^{\ast}\hat{H}(t))\\
&  +\hat{q}^{0}(t)\tilde{\sigma}(t)+%
{\textstyle\int_{\mathbb{R}_{0}}}
\hat{r}^{0}(t,\zeta)\tilde{\gamma}(t,\zeta)\nu(d\zeta)\}dt].
\end{align*}
Similarly, we obtain
\begin{align}
&  \mathbb{E}[\langle\hat{p}_{1}^{1}(T),\tilde{M}(T)\rangle]\nonumber\\
&  =\mathbb{E}[%
{\textstyle\int_{0}^{T}}
\langle\hat{p}_{1}^{1}(t),d\tilde{M}(t)\rangle+%
{\textstyle\int_{0}^{T}}
\tilde{M}(t)d\tilde{p}_{1}^{1}(t)]\nonumber\\
&  =\mathbb{E}[%
{\textstyle\int_{0}^{T}}
\langle\hat{p}_{1}^{1}(t),\tilde{M}^{\prime}(t)\rangle dt-{%
{\textstyle\int_{0}^{T}}
\{}\langle\nabla_{m}\hat{H}_{1}(t),\tilde{M}(t)\rangle-\nabla_{\bar{m}}^{\ast
}\hat{H}_{1}(t)\tilde{M}(t)\}dt].\nonumber
\end{align}
In the above we have used that the expectation of the martingale terms, i.e.
the $dB(t)$- and $\tilde{N}(dt,d\zeta)$ -integrals, have mean zero. Taking the
limit superior \color{black} and using the transversality
conditions $\left(  \ref{tcond_1}\right)  $ combined with \eqref{eq3.13} we
obtain, using that $u$ and $\hat{u}$ are $\mathbb{G}$-adapted,
\begin{align*}
A  &  \leq-\underset{T\rightarrow\infty}{\overline{\lim}}\mathbb{E}[\hat
{p}^{0}(T)\tilde{X}(T)]-\underset{T\rightarrow\infty}{\overline{\lim}%
}\mathbb{E}[\hat{p}^{1}(T)\tilde{M}(T)]+\mathbb{E[}{{%
{\textstyle\int_{0}^{T}}
}}\tfrac{\partial\hat{H}}{\partial u}(t)_{u=\hat{u}(t)}\tilde{u}(t)dt]\\
&  \leq\mathbb{E[}{{%
{\textstyle\int_{0}^{T}}
}}\tfrac{\partial\hat{H}}{\partial u}(t)_{u=\hat{u}(t)}\tilde{u}(t)dt]\leq0,
\end{align*}
since $u\mapsto\mathbb{E}[\hat{H}(t,u)|\mathcal{G}_{t}]$ is maximal at
$u=\hat{u}(t)$. That completes the proof. \hfill$\square$ \bigskip

\subsection{Necessary maximum principle under partial information}

\noindent We now consider the converse, i.e. we look for necessary conditions
of optimality. The following result is the infinite horizon version of Theorem
\ref{nece-fin}:

\begin{theorem}
Assume that Assumptions A1-A3 of Section 3.2 hold but now with $t\in
\lbrack0,\infty)$. Let $u\in\mathcal{U}_{adm}$ with corresponding solutions
$X$ and $(p^{0},q^{0},r^{0})$ and $(p^{1},q^{1},r^{1})$ of the forward and the
backward stochastic differential equations $\left(  \ref{F}\right)  $ and
\eqref{eq4.11} and \eqref{eq4.11'} respectively with corresponding derivative
process $Z$ given by \eqref{eq3.25} but now with the time horizon $[0,\infty)$.

Moreover, assume that the following transversality condition holds:%
\begin{equation}
\underset{T\rightarrow\infty}{\overline{\lim}}\mathbb{E}[p^{0}%
(T)Z(T)]=\underset{T\rightarrow\infty}{\overline{\lim}}\mathbb{E}[\left\langle
p^{1}(T),DM(T)\right\rangle ]=0;\quad\text{ for all bounded }\pi\in
\mathcal{U}_{adm}. \label{trv_c_n}%
\end{equation}

Then the following are equivalent:

\begin{itemize}
\item
\begin{equation}
\label{eq4.18}\tfrac{d}{d\lambda}J(u+\lambda\pi)|_{\lambda=0}=0 \text{ for all
bounded } \pi\in\mathcal{U}_{adm}.
\end{equation}

\item
\begin{equation}
\mathbb{E[}\tfrac{\partial H}{\partial u}(t,u)|\mathcal{G}_{t}]=0.
\label{eq4.19}%
\end{equation}

\end{itemize}
\end{theorem}

\noindent{Proof.} \quad Assume that \eqref{eq4.18} holds. Then%
\begin{align}
0  &  =\tfrac{d}{d\lambda}J(u+\lambda\pi)|_{\lambda=0}\label{eq4.12}\\
&  =\mathbb{E[}%
{\textstyle\int_{0}^{\infty}}
\{\tfrac{\partial f}{\partial x}(t)Z(t)+\langle\nabla_{\overline{x}}%
f(t),Z_{t}\rangle+\langle\nabla_{m}f(t),DM(t)\rangle\nonumber\\
&  +\langle\nabla_{\overline{m}}f(t),DM_{t}\rangle+\tfrac{\partial f}{\partial
u}(t)\pi(t)\}dt].\nonumber
\end{align}
By the definition of the Hamiltonian \eqref{eq3.6}, we have%
\[%
\begin{array}
[c]{cc}%
\nabla f(t) & =\nabla H(t)-\nabla b(t)p^{0}(t)-\nabla\sigma(t)q^{0}(t)-%
{\textstyle\int_{\mathbb{R}_{0}}}
\nabla\gamma(t,\zeta)r^{0}(t,\zeta)\nu(d\zeta),
\end{array}
\]
where $\nabla=(\tfrac{\partial}{\partial x},\nabla_{\overline{x}},\nabla
_{m},\nabla_{\overline{m}},\frac{\partial}{\partial u}).$

\noindent Define a sequence of stopping times by%
\[%
\begin{array}
[c]{c}%
\tau_{n}(\cdot):=T\wedge\inf\{t\geq0:\int_{0}^{t}\left\{  (p^{0}%
(s))^{2}+(q^{0}(s))^{2}\right. \\
+%
{\textstyle\int_{\mathbb{R}_{0}}}
(r^{0}(s,\zeta))^{2}\nu(d\zeta)+\pi^{2}(s))ds\geq n\},n\in%
\mathbb{N}
.
\end{array}
\]
Clearly $\tau_{n}\rightarrow T$ $\mathbb{P}$-a.s. as $n\rightarrow\infty$.
Applying the Itô formula, we get%
\[%
\begin{array}
[c]{l}%
\mathbb{E}[p^{0}({\normalsize T})Z({\normalsize T})]+\mathbb{E}[\left\langle
p^{1}(T),DM(T)\right\rangle ]=\underset{n\rightarrow\infty}{\lim}%
(\mathbb{E}[p^{0}({\normalsize \tau_{n}})Z({\normalsize \tau_{n}}%
)]+\mathbb{E}[\left\langle p^{1}({\normalsize \tau_{n}}),DM({\normalsize \tau
_{n}})\right\rangle ])\\
=\mathbb{E}\Big[%
{\textstyle\int_{0}^{{\normalsize \tau_{n}}}}
p^{0}(t)\Big\{\tfrac{\partial b}{\partial x}(t)Z(t)+\langle\nabla
_{\overline{x}}b(t),Z_{t}\rangle+\langle\nabla_{m}b(t),DM(t)\rangle
+\langle\nabla_{\overline{m}}b(t),DM_{t}\rangle+\tfrac{\partial b}{\partial
u}(t)\pi(t)\\
-\big(\tfrac{\partial H}{\partial x}(t)+\mathbb{E}[\nabla_{\bar{x}}^{\ast
}H(t)|\mathcal{F}_{t}]\big)Z(t)\\
+q(t)\big(\tfrac{\partial\sigma}{\partial x}(t)Z(t)+\langle\nabla
_{\overline{x}}\sigma(t),Z_{t}\rangle+\langle\nabla_{m}\sigma(t),DM(t)\rangle
+\langle\nabla_{\overline{m}}\sigma(t),DM_{t}\rangle+\tfrac{\partial\sigma
}{\partial u}(t)\pi(t)\big)\\
+%
{\textstyle\int_{\mathbb{R}_{0}}}
r(t,\zeta)\big(\tfrac{\partial\gamma}{\partial x}(t,\zeta)Z(t)+\langle
\nabla_{\overline{x}}\gamma(t,\zeta),Z_{t}\rangle+\langle\nabla_{m}%
\gamma(t,\zeta),DM(t)\rangle\\
+\langle\nabla_{\overline{m}}\gamma(t,\zeta),DM_{t}\rangle+\tfrac
{\partial\gamma}{\partial u}(t,\zeta)\pi(t)\big)\nu(d\zeta)\Big\}dt\Big]\\
+\mathbb{E}\Big[%
{\textstyle\int_{0}^{{\normalsize \tau_{n}}}}
\{\left\langle p^{1}(t),DM^{\prime}((t)\right\rangle -\left\langle \nabla
_{m}H(t),DM((t)\right\rangle -\nabla_{\bar{m}}^{\ast}H(t)DM((t)\}dt].
\end{array}
\]
Taking the limit superior, \color{black} combining this with
\eqref{eq4.12} and using the transversality condition $\left(  \ref{trv_c_n}%
\right)  $, we get
\[
0=\underset{T\rightarrow\infty}{\overline{\lim}}\mathbb{E}[p^{0}%
({\normalsize T})Z({\normalsize T})]+\underset{T\rightarrow\infty}%
{\overline{\lim}}\mathbb{E}[\left\langle p^{1}(T),DM(T)\right\rangle
]=\mathbb{E[}%
{\textstyle\int_{0}^{{\normalsize \infty}}}
\tfrac{\partial H}{\partial u}(t)\pi(t)dt]\text{.}%
\]
\newline Now choose $\pi(t)=\alpha\mathbf{1}_{\left(  t_{0},T\right]  }(t)$,
where $\alpha=\alpha(\omega)$ is bounded and $\mathcal{G}_{t_{0}}$-measurable
and $t_{0}\in\lbrack0,T)$. Then we deduce that
\[%
\begin{array}
[c]{ll}%
\mathbb{E[}%
{\textstyle\int_{t_{0}}^{\infty}}
\tfrac{\partial H}{\partial u}(t)\alpha dt] & =0\text{.}%
\end{array}
\]
Differentiating with respect to $t_{0}$ we obtain%
\[%
\begin{array}
[c]{ll}%
\mathbb{E[}\tfrac{\partial H}{\partial u}(t_{0})\alpha] & =0\text{.}%
\end{array}
\]
Since this holds for all such $\alpha$, we conclude that%
\[%
\begin{array}
[c]{ll}%
\mathbb{E[}\tfrac{\partial H}{\partial u}(t_{0})|\mathcal{G}_{t_{0}%
}]=0\text{,} & \text{ which is \eqref{eq4.19}.}%
\end{array}
\]
This argument can be reversed, to prove that \eqref{eq4.19} $\Longrightarrow$
\eqref{eq4.18}. We omit the details.$\qquad\square$ \newline

\section{Optimal consumption from an elephant memory cash flow}

\noindent To illustrate our results, let us consider an example of an infinite
horizon optimal consumption problem, where the wealth process of the investor
$X=(X^{u}(t))_{t\geq0}$ is given by the following dynamics:%

\begin{equation}
\left\{
\begin{array}
[c]{ll}%
dX^{u}(t) & =\big\{\left\langle F,X_{t}^{u}\right\rangle -u(t)\big\}dt+\beta
X^{u}(t)dB(t);\,t\geq0,\\
X^{u}(0) & =x_{0}>0,
\end{array}
\right.  \nonumber\label{eq5.1}%
\end{equation}
where $u(t)\geq0$ denotes the consumption rate (our control), $\beta>0$
(constant) denotes the volatility and $F(\cdot):L_{0}(\mathbb{R}%
)\mapsto\mathbb{R}$ is a bounded linear operator on the whole memory path
$X_{t}^{u}=\{X^{u}(t-s)\}_{0\leq s\leq t}$ of $X$ up to time $t$. 
Thus the term $\left\langle F,X_{t}^{u}\right\rangle $ represents a drift term
in the dynamics depending on the whole history of the process. A specific
example is given below. \color{black} \newline We define $\mathcal{U}_{adm}$
to be the set of nonnegative adapted processes $u$ such that%

\[
\mathbb{E}\left[  {%
{\textstyle\int_{0}^{\infty}}
|X^{u}(t)|}^{2}dt\right]  <\infty
\]
For $u\in\mathcal{U}_{adm}$ we also require that $u$ satisfies the following
\emph{budget constraint}:\newline The expected total discounted consumption is
bounded by the initial capital $x_{0}$, i.e.:%

\begin{equation}
\mathbb{E}\left[  \int_{0}^{\infty}e^{-\rho t}u(t)dt\right]  \leq x_{0},
\label{eq5.2}%
\end{equation}
where $\rho>0$ is a given discount exponent. Consider the following problem:

\begin{problem}
Find $\hat{u}\in\mathcal{U}_{adm}$ such that
\begin{equation}
J(\hat{u})=\sup_{u\in\mathcal{U}_{adm}}J(u),
\end{equation}
where the performance functional $J(u)$ is the total discounted logarithmic
utility of the consumption $u$, i.e.
\begin{equation}
J(u)=\mathbb{E}\left[  {{\int_{0}^{\infty}}}e^{-\delta t}\ln(u(t))dt\right]
;u\in\mathcal{U}_{adm},
\end{equation}
for some constant $\delta>0.$
\end{problem}

\noindent The Hamiltonian in this case takes the form%
\[
H(t,x,\bar{x},u,p^{0},q^{0})=e^{-\delta t}\ln(u)+p^{0}[\left\langle F,\bar
{x}\right\rangle -u]+q^{0}\beta x,
\]
and the adjoint process pair $\left(  p^{0}(t),q^{0}(t)\right)  $ is a
solution of the corresponding adjoint BSDE%

\begin{equation}
dp^{0}(t)=-\big\{\beta q^{0}(t)+\mathbb{E}[\nabla_{\bar{x}}^{\ast
}H(t)|\mathcal{F}_{t}]\big\}dt+q^{0}(t)dB(t);t\in\lbrack0,\infty).
\label{eq5.5}%
\end{equation}

\noindent Note that by Lemma \ref{lem} we have
\[
\nabla_{\bar{x}}^{\ast}H(t)=\left\langle F,(p^{0})^{t}\right\rangle .
\]
For example, let us from now on assume that $F(\cdot)$ is a weighted average
operator of the form
\begin{equation}
\left\langle F,\bar{x}\right\rangle ={%
{\textstyle\int_{0}^{t}}
}e^{-\rho r}x(r)dr. \label{average}%
\end{equation}
Then we get
\[
\nabla_{\bar{x}}^{\ast}H(t)={%
{\textstyle\int_{0}^{\infty}}
}e^{-\rho r}p^{0}(t+r)dr,
\]
and the state equation becomes
\begin{equation}
\left\{
\begin{array}
[c]{ll}%
dX^{u}(t) & =\big\{\int_{0}^{t}e^{-\rho r}X^{u}(t-r)dr-u(t)\big\}dt+\beta
X^{u}(t)dB(t);\,t\geq0,\\
X^{u}(0) & =x_{0}>0.
\end{array}
\right.  \nonumber\label{eq5.3a}%
\end{equation}
The adjoint BSDE \eqref{eq5.1} will take the form
\begin{equation}
dp^{0}(t)=-\big\{\beta q^{0}(t)+\mathbb{E}[{%
{\textstyle\int_{0}^{\infty}}
}e^{-\rho r}p^{0}(t+r)dr|\mathcal{F}_{t}]\big\}dt+q^{0}(t)dB(t);t\in
\ [0,\infty). \label{eq5.2}%
\end{equation}
Maximising the Hamiltonian with respect to $u$ gives the following equation
for a possible optimal consumption rate $u=\hat{u}$:
\[
e^{-\delta t}\tfrac{1}{\hat{u}(t)}-\hat{p}^{0}(t)=0,
\]
i.e.
\begin{equation}
\hat{u}(t)=\tfrac{e^{-\delta t}}{\hat{p}^{0}(t)}. \label{eq5.3}%
\end{equation}
With this choice $u=\hat{u}$ the equations above get the form
\[
\left\{
\begin{array}
[c]{ll}%
d\hat{X}(t) & =\big\{\int_{0}^{t}e^{-\rho r}\hat{X}(t-r)dr-\frac{e^{-\delta
t}}{\hat{p}(t)}\big\}dt+\beta\hat{X}(t)dB(t);\,t\geq0,\\
\hat{X}(0) & =x_{0},
\end{array}
\right.
\]
and
\begin{equation}
d\hat{p}^{0}(t)=-\big\{\beta\hat{q}^{0}(t)+\mathbb{E}[{{\textstyle\int
_{0}^{\infty}}}e^{-\rho r}\hat{p}^{0}(t+r)dr|\mathcal{F}_{t}]\big\}dt+\hat
{q}^{0}(t)dB(t);t\in\ [0,\infty), \label{eq5.2}%
\end{equation}
We want to find a solution $(\hat{p}^{0},\hat{q}^{0})$ of this infinite
horizon BSDE such that the transversality condition holds, i.e.
\[
\overline{\lim_{T\rightarrow\infty}}\mathbb{E}[\hat{p}^{0}(T)(X^{u}(T)-\hat
{X}(T))]\geq0
\]
for all admissible controls $u$.\newline

\begin{remark}
This problem may be regarded as an infinite horizon version of a stochastic
control problem for a Volterra equation, without memory. To see this, note
that by a change of variable and a change of the order of integration the
equation \eqref{eq5.3a} can be written
\begin{align}
X(t) &  =x_{0}+\int_{0}^{t}\left(  \int_{0}^{s}e^{-\rho r}X(t-r)dr\right)
ds-\int_{0}^{t}u(s)ds+\int_{0}^{t}\beta X(s)dB(s)\nonumber\label{eq5.8}\\
&  =x_{0}+\int_{0}^{t}\frac{1}{\rho}(1-e^{\rho(s-t)})X(s)ds-\int_{0}%
^{t}u(s)ds+\int_{0}^{t}\beta X(s)dB(s),
\end{align}
which is a stochastic Volterra equation of the type studied in \cite{AO1} and
\cite{AOY}.
\end{remark}

Let us try to assume that $\hat{q}^{0}=0$ and hence that $\hat{p}^{0}$ is
deterministic. Then the equation for $\hat{p}^{0}(t)$ reduces to the integral
equation
\begin{equation}
d\hat{p}^{0}(t)=-\Big(\int_{0}^{\infty}e^{-\rho r}\hat{p}^{0}(t+r)dr\Big)dt.
\end{equation}
 By a similar procedure as in \eqref{eq5.8} above we obtain that
this equation can be transformed to the equation \color{black}
\begin{equation}
d\hat{p}^{0}(t)=-\frac{1}{\rho}(1-e^{-\rho t})\hat{p}^{0}(t)dt,
\end{equation}
which has the general solution
\begin{equation}
\hat{p}^{0}(t)=\hat{p}^{0}(0)\exp(-\frac{t}{\rho}+\frac{1-e^{-\rho t}}%
{\rho^{2}});\quad t\geq0, \label{eq5.11}%
\end{equation}
for some constant $\hat{p}^{0}(0)$.

\noindent Substituted into \eqref{eq5.3a} this gives
\begin{equation}
\hat{u}(t)=\frac{1}{\hat{p}^{0}(0)}\exp\Big((\frac{1}{\rho}-\delta
)t-\frac{1-e^{-\rho t}}{\rho^{2}}\Big). \label{eq5.12}%
\end{equation}

\noindent The problem is to find $\hat{p}^{0}(0)$ such that the following two
conditions hold:%

\begin{align}
\overline{\lim_{T\rightarrow\infty}}\hat{p}^{0}(T)\mathbb{E}[\hat{X}(T)] &
=0\label{eq5.13}\\
\mathbb{E}\left[  \int_{0}^{\infty}e^{-\rho t}u(t)dt\right]   &  \leq
x_{0}\quad\text{(the budget constraint).}\label{eq5.14}%
\end{align}
Define $y_{0}(t)$ to be the solution of the integral equation
\begin{equation}
y_{0}(t)=x_{0}+\int_{0}^{t}\frac{1}{\rho}(1-e^{\rho(s-t)})y_{0}(s)ds;\quad
t\geq0,
\end{equation}
and let $\lambda_{0}>0$ be the top Lyapunov exponent of $y_{0}.$ See e.g.
Kunita \cite{K} and Mang and Sheng \cite{MS} for more information about
Lyapunov exponents. Then, since clearly
\[
y_{0}(t)\geq\mathbb{E}[\hat{X}(t)]\text{ for all }t\geq0,
\]
we see by \eqref{eq5.11} that if
\begin{equation}
\rho<\frac{1}{\lambda_{0}},
\end{equation}
then
\begin{equation}
\overline{\lim_{T\mapsto\infty}}\hat{p}^{0}(T)\mathbb{E}[\hat{X}(T)]=0,
\end{equation}
and hence \eqref{eq5.13} holds for any choice of $\hat{p}^{0}(0)$. By
\eqref{eq5.12} the budget constraint \eqref{eq5.14} gives
\begin{equation}
\hat{p}^{0}(0)\geq\frac{1}{x_{0}}\int_{0}^{\infty}\exp\Big((\frac{1}{\rho
}-\delta-\rho)t-\frac{1-e^{-\rho t}}{\rho^{2}}\Big)dt.\label{eq5.22}%
\end{equation}
The admissible value of $\hat{p}^{0}(0)$ that gives the maximal consumption is
therefore, by \eqref{eq5.12},
\begin{equation}
\hat{p}^{0}(0)=\frac{1}{x_{0}}\int_{0}^{\infty}\exp\Big((\frac{1}{\rho}%
-\delta-\rho)t-\frac{1-e^{-\rho t}}{\rho^{2}}\Big)dt.\label{eq5.23}%
\end{equation}
We summarise what we have proved as follows:

\begin{theorem}
Assume that
\[
\rho<\frac{1}{\lambda_{0}}.
\]
Then the optimal consumption rate $\hat{u}(t)$ for Problem 5.1, with $F$
defined by \eqref{average}, is given by \eqref{eq5.12}, where $\hat{p}^{0}(0)$
is given by \eqref{eq5.23}.
\end{theorem}

\end{document}